\documentclass[psamsfonts,12pt]{article}

\usepackage{xspace}
\usepackage[ps,tips,matrix,arrow]{xy}
\usepackage{amsmath,amssymb,amsxtra,amsthm}

\newtheorem{theorem}{Theorem}[section]

\newtheorem{proposition}[theorem]{Proposition}

\newtheorem{definition}[theorem]{Definition}

\newtheorem{remark}[theorem]{Remark}

\newcommand{\cat}[1]{\ensuremath{\mathsf{#1}}}
\newcommand{\Z}{\ensuremath{\mathbb{Z}}}

\begin{document}

\title{A simplicial model for the Hopf map}

\author{Orin R. Sauvageot\footnote{\'Ecole Polytechnique F\'ed\'erale de Lausanne --
 Institute of Geometry, Algebra and Topology -- CH-1015~Lausanne, Switzerland}\\
\texttt{orin.sauvageot@epfl.ch}}

\maketitle

\begin{abstract}

We give an explicit simplicial model for the Hopf map $S^3\to S^2$. For this
purpose, we construct a model of $S^3$ as a principal twisted cartesian
product
$K\times_{\eta}S^2$, where
$K$ is a simplicial model for $S^1$ acting by left multiplication on itself,
$S^2$ is given the simplest simplicial model and the twisting map is 
$\eta:\left(S^2\right)_n \to (K)_{n-1}$. We construct a Kan complex for the
simplicial model
$K$ of
$S^1$. The simplicial model for the Hopf map is then the projection
$K\times_{\eta}S^2\to S^2$.

\end{abstract}

\noindent\textit{Mathematics Subject Classification~:} 55U10, 55Q40; 55R10 

\section{Introduction}

The motivation for finding a simplicial model for the Hopf map arose when trying to find a simple test to decide
whether the stabilisation of a certain interesting model category $\cat{M}$ is different from that of the category of
chain complexes of abelian groups. As detailed in \cite[chapter 6]{SauPhD}, consider the following situation. Let
$\cat{M}$ be a symmetric monoidal model category whose stabilisation exists and suppose there is a monoidal Quillen
adjunction
$F:\cat{sS}\leftrightarrow \cat{M}:G$ between the category of simplicial sets and $\cat{M}$. In the stable category of
chain complexes, the Hopf map vanishes. Therefore, if we have a good simplicial model for the Hopf map that allows us
to show that the multiply suspended images under the functor $F$ of this simplicial model never vanish, then the
stabilisation of $\cat{M}$ is different from that of chain complexes.

The main result of this paper is that a very good simplicial model of the Hopf map is the projection
$p:K\times_{\eta} S^2\to S^2$ of a principal twisted cartesian product of a simplicial model $K$ of $S^1$ with
the simplest simplicial model for $S^2$. The proof of this result shows that we are able to model simplicially any
$S^1$-bundle of base $S^2$.

This paper is structured in the following manner. Section \ref{S:PTCP} recalls the notions and results related to
principal twisted cartesian products. In section \ref{S:model-S1} we construct a Kan model for $S^1$, which is required
to carry enough structure. Section \ref{S:computations} gives explicit computations of the Kan model as well as of the
twisting map. Finally, we prove the main result in section \ref{S:main}.

We are indebted to Kathryn Hess and Andrew Tonks for the idea of using principal twisted cartesian products. The idea
of using the Hopf map as a test came out from a discussion with Stefan Schwede. This work has been carried out with
the financial support of the Swiss National Science Foundation.

\section{Principal twisted cartesian products}\label{S:PTCP} 
This section is devoted to explaining the tools for building a
simplicial model for $S^3$ with enough structure to capture the Hopf map.
\begin{definition}
Let $F$ and $B$ be two simplicial sets. Let $H$ be a simplicial group acting on
the left on $F$.   Let $\zeta:B\to H$ be a map of graded sets of degree $-1$ such that $\zeta_n	:B_n
\to H_{n-1}$ satisfies the following identities: 
\begin{eqnarray*}
		 \partial_{0}\zeta(b) &=& 
		 \left(\zeta(\partial_{0}b)\right)^{-1}\zeta(\partial_{1}b) \\
		 \partial_{i}\zeta(b) &=& \zeta(\partial_{i+1}b) \quad\text{for } i>0 \\
		 s_{i}\zeta(b) &=& \zeta(s_{i+1}b)\quad\text{for } i\geq 0 \\
		 \zeta(s_{0}b) &=& id_{n}\quad\text{for } b\in B_{n}.
	\end{eqnarray*}

The map
$\zeta$ is the
\textsl{twisting map}. A
\textsl{twisted cartesian product} of fibre
$F$, base
$B$ and group $H$ is a simplicial set denoted $F\times_\zeta B$ satisfying
$$
(F\times_\zeta B)_n = F_n	\times B_n	
$$
with faces and degeneracies as follows~:
\begin{enumerate}
\item $\partial_i	(f,b)=(\partial_i f, \partial_i b)$ for $i>0$
\item $\partial_0 (f,b))= (\zeta(b)\partial_0 f,\partial_0 b)$
\item $s_i(f,b)= (s_i f,s_i	b	)$ for $i\geq 0$.
\end{enumerate}

\noindent Furthermore, if $F=H$ acting on itself by left multiplication, then $F\times_\zeta B$ is a \textsl{principal
twisted cartesian product (PTCP)}. 
\end{definition}
We will also use the terminology
``twisted cartesian product'' for the projection $p:F\times_\zeta B\to B$.


The following proposition is a classical result whose proof can be found in
\cite[proposition 18.4]{MaySS}.

\begin{proposition}\label{P:fib-princ}
Let $p:F\times_\zeta B\to B$ be a twisted cartesian product with group $H$.
If the fiber $F$ is a Kan complex, then~:
\begin{enumerate}
\item the projection $p$ is a Kan fibration, and
\item if $F=H$, $p$ is a principal fibration.
\end{enumerate}
\end{proposition}

\begin{remark}\label{R:identities}
Let $S:\cat{TOP}\to \cat{sS}$ be the singular functor from the category
$\cat{TOP}$ of topological spaces to the category $\cat{sS}$ of simplicial sets.
A map $f$ is a (Serre) fibration if and only if $S(f)$ is a Kan fibration. Thus,
a principal fibration (or fibre bundle) in $\cat{TOP}$ passes via the functor $S$
to a principal fibration in $\cat{sS}$. Since the Hopf map $S^3\to S^2$ is a
fibration in $\cat{TOP}$, the corresponding simplicial model has to be a Kan
fibration. As a consequence, if we want to model $S^3$ as a PTCP $K\times_\eta S^2\to S^2$, this has to be a Kan
fibration, which it is when $K$ is a Kan complex, by proposition \ref{P:fib-princ}. 
We construct such a PTCP in the following sections.
\end{remark}

\section{The simplicial model for $S^1$}\label{S:model-S1}
In short, to build a Kan model of $S^1$ we let $\Z(2)$ denote a chain
complex concentrated in degree two, and we apply a functor $\Gamma$ to obtain a
simplicial abelian group $\Gamma\Z(2)$. By applying the loop
group functor $G$, the model of
$S^1$ is given by $G\Gamma\Z(2)$. The
latter is always a Kan complex, since every simplicial group is a Kan complex. More precisely we give the following
definitions.

\begin{definition}
Let $\cat{sAb}$ be the category of simplicial abelian groups and let $\cat{CC}$ be the
category of chain complexes of abelian groups. We define the functor
$\Gamma:\cat{CC}\to\cat{sAb}$ as follows. For any $(X,\partial)\in\cat{CC}$, the
simplicial abelian group $\Gamma(X)$ is given by~:

\begin{enumerate}
		\item \begin{equation}
\Gamma_{n}(X)=
X_{n}\bigoplus_{r=0}^{n-1}\sum_{k=n-r}\sigma_{j_{k}}\ldots\sigma_{j_{1}} 
		X_{r}
\end{equation}\label{E:Gamma} where $\sigma_{j_{k}}\ldots\sigma_{j_{1}}X_{r}$ is
the 
		abelian group whose elements are the symbols
		$\sigma_{j_{k}}\ldots\sigma_{j_{1}}x$ with $x\in X_{r}$. 
		The sum $\sum_{k=n-r}$ is taken over all sequences of indices 
		$\{j_{i}\}$ such that $0\leq j_{1}<j_{2}<\cdots<j_{k}<n$. 
		
		The addition of symbols is defined by
		$$
		\sigma_{j_{k}}\ldots\sigma_{j_{1}}x +\sigma_{j_{k}}\ldots\sigma_{j_{1}}y = \sigma_{j_{k}}\ldots\sigma_{j_{1}}(x+y).
		$$
		Degeneracies and faces are given by~:

\vspace{.2cm}
		
\item $s_{i}:\Gamma_{n}(X)\to \Gamma_{n+1}(X)$ is defined by
	  \begin{enumerate}
		  \item $s_{i}x=\sigma_{i}x$ for $x\in X_{n}$
		  \item if $k=n-r$ and $x\in X_{r}$ then
		  $$
		  s_{i}\sigma_{j_{k}}\ldots\sigma_{j_{1}}x=\sigma_{h_{k+1}}\ldots\sigma_{h_{1}}x
		  $$
		 when $s_{i}s_{j_{k}}\ldots s_{j_{1}}=s_{h_{k+1}}\ldots 
		  s_{h_{1}}$ and where $s_{h_{k+1}}\ldots s_{h_{1}}$ is written in 
		  the canonical form\footnote{Every composition of degeneracies and/or faces can
be written in the canonical form with the aid of the simplicial identities.}, i.e.
$h_{k+1}>h_k>...>h_1$.
	  \end{enumerate}

  \item $\partial_{i}:\Gamma_{n}(X)\to\Gamma_{n-1}(X)$ is defined by
		\begin{enumerate}
			\item $\partial_{n}x=\partial(x)$ and $\partial_{i}x=0$ if $i<n$ 
			and $x\in X_{n}$.
			\item if $k=n-r$ and $x\in X_{r}$ then
			$$
			\hspace{-1.5cm}\partial_{i}\sigma_{j_{k}}\ldots\sigma_{j_{1}}x=\begin{cases}
			\sigma_{h_{k-1}}\ldots\sigma_{h_{1}}x \\
			\sigma_{h_{k}}\ldots\sigma_{h_{1}}\partial(x)\\
			0 \end{cases}$$

			if respectively 
$$
\partial_{i}s_{j_{k}}\ldots 
			s_{j_{1}}=\begin{cases}
			s_{h_{k-1}}\ldots s_{h_{1}}\\
			s_{h_{k}}\ldots s_{h_{1}}\partial_{r}\\
			s_{h_{k}}\ldots s_{h_{1}}\partial_{j}\ \ j<r
			\end{cases}
	        $$
			where the right hand side is written in the 
			canonical form.
		\end{enumerate}	
		
	\end{enumerate}
\end{definition}

\noindent We now define the functor $G$.

\begin{definition}\label{D:G}
Let $\cat{sGr}$ the category of simplicial groups and let $K$ be a simplicial
set. We define the functor
$G:\cat{sS}\to
\cat{sGr}$ as follows. The group $G_n (K)=G(K)_n$ is the free group generated by
the elements of $K_{n+1}$ modulo the relations $s_0 x=id_n$ for all $x\in K_n$.

If $x\in K_{n+1}$, let $\zeta(x)$ be the class of $x$ in $G_n(K)$. Faces and
degeneracies of $G(K)$ are defined on generators by the relations~:
\begin{eqnarray}\label{E:FD}
\zeta(\partial_0 x)\partial_0 \zeta(x) &=& \zeta(\partial_1 x)\label{E:FD1} \\
\partial_i \zeta(x) &=& \zeta(\partial_{i+1}x)\qquad\text{if }i>0\label{E:FD2}\\
s_i	\zeta(x) &=& \zeta(s_{i+1}x)\qquad\text{if }i\geq 0.\label{E:FD3}
\end{eqnarray}
\end{definition}

By extension we have homomorphisms $\partial_i:G_n(K)\to G_{n-1}(K)$ and
$s_i:G_n(K)\to G_{n+1}(K)$. Clearly, $G(K)$ is a simplicial group.

\begin{remark}\label{R:PTCP}
The morphism $\zeta$ of definition \ref{D:G} is clearly a twisting map. 
Hence, for every simplicial abelian group
$K$ we have a twisted cartesian product $G(K)\times_{\zeta}K$, which is acyclic. The reader may refer to
\cite[pp. 118--123]{MaySS} for details.
\end{remark}

By \cite[Remarks 23.7]{MaySS}, $\Gamma\Z(2)$ is a $K(\Z,2)$, hence a simplicial model for
$BS^1$. $G\Gamma\Z(2)$ is then a model for $\Omega BS^1$, hence for $S^1$.

\section{Some computations}\label{S:computations}
This section is devoted to clarifying the previous construction by giving explicit
computations of $G\Gamma \Z(2)$ and the map $\eta$. For this we will choose a
simplicial model for $S^2$ consisting in one non degenerate simplex in degree two and only
degeneracies above.

To compute $\Gamma_n \Z(2)$ we use formula (\ref{E:Gamma}). Since $\Z(2)$ is
concentrated in degree two, we obtain 
\begin{equation}\label{E:Gamma-Z-2}
\Gamma_n \Z(2)= \bigoplus_{0\leq j_1<\cdots<j_{n-2}<n}
\sigma_{j_{n-2}}\cdots\sigma_{j_1} \Z.
\end{equation}
As an example, in degree three, the faces and degeneracies are given for all
$z\in \Z$ by
\begin{alignat*}{4}
 \partial_{0}(\sigma_{0}z) &= z \quad \partial_{1}(\sigma_{0}z) &= z 
 \quad \partial_{2}(\sigma_{0}z) &= 0 \quad \partial_{3}(\sigma_{0}z) &= 0 \\ 
 \partial_{0}(\sigma_{1}z) &= 0 \quad \partial_{1}(\sigma_{1}z) &= z 
 \quad \partial_{2}(\sigma_{1}z) &= z \quad \partial_{3}(\sigma_{1}z) &= 0 \\ 
 \partial_{0}(\sigma_{2}z) &= 0 \quad \partial_{1}(\sigma_{2}z) &= 0 
 \quad \partial_{2}(\sigma_{2}z) &= z \quad \partial_{3}(\sigma_{2}z) &= z 
\end{alignat*} 

\begin{alignat*}{4}
 s_{0}(\sigma_{0}z) &= \sigma_{1}\sigma_{0}z  \quad 
 s_{1}(\sigma_{0}z) &=  \sigma_{1}\sigma_{0}z 
 \quad s_{2}(\sigma_{0}z) &= \sigma_{2}\sigma_{0}z  \quad 
 s_{3}(\sigma_{0}z) &= \sigma_{3}\sigma_{0}z  \\ 
 s_{0}(\sigma_{1}z) &= \sigma_{2}\sigma_{0}z   \quad s_{1}(\sigma_{1}z) &= 
 \sigma_{2}\sigma_{1}z   
 \quad s_{2}(\sigma_{1}z) &= \sigma_{2}\sigma_{1}z  \quad 
 s_{3}(\sigma_{1}z) &= \sigma_{3}\sigma_{1}z \\ 
 s_{0}(\sigma_{2}z) &= \sigma_{3}\sigma_{0}z  \quad 
 s_{1}(\sigma_{2}z) &= \sigma_{3}\sigma_{1}z  
 \quad s_{2}(\sigma_{2}z) &= \sigma_{3}\sigma_{2}z  \quad 
 s_{3}(\sigma_{2}z) &= \sigma_{3}\sigma_{2}z.
\end{alignat*}

\noindent For $G\Gamma \Z(2)$, we have
\begin{eqnarray}
(G\Gamma \Z(2))_{0} &=& \{e\},\qquad (G\Gamma \Z(2))_{1}=\mathcal{F}\{
\Z\backslash
\{0\}\}\notag\\ 
(G\Gamma \Z(2))_{2} &=& \mathcal{F}\{\sigma_{2} \Z\oplus\sigma_{1} \Z\}\notag\\
(G\Gamma \Z(2))_{3}
&=&\mathcal{F}\{\sigma_{2}\sigma_{1}
\Z\oplus\sigma_{3}\sigma_{2} \Z\oplus\sigma_{3}\sigma_{1} \Z\}\notag\\ 
(G\Gamma \Z(2))_{4}
&=&\mathcal{F}\{\sigma_{4}\sigma_{2}\sigma_{1} \Z\oplus
\sigma_{4}\sigma_{3}\sigma_{2}
\Z\oplus\sigma_{4}\sigma_{3}\sigma_{1} \Z\oplus\sigma_{3}\sigma_{2}\sigma_{1}
\Z\}\notag\\
 &\vdots&\notag\\
 (G\Gamma \Z(2))_{n} &=&
\mathcal{F}\left\{\bigoplus_{0<j_{1}<\cdots<j_{n-1}<n+1}\sigma_{j_{n-1}}
\ldots\sigma_{j_{1}} \Z\right\}.\label{E:G-Gamma}
\end{eqnarray}	
where $\mathcal{F}\{\}$ stands for the free group generated by elements inside
$\{\}$. Notice that $s_{0}(\sigma_{j_{n-1}}\ldots\sigma_{j_{1}}z)$ 
can always be expressed in a form ending by $\sigma_{0}z$. Hence each term
containing
$\sigma_{0} \Z$ is trivial and gives the first strict inequality in
$0<j_{i}<\cdots<j_{n-1}<n+1$. Faces and degeneracies are given by the formulae
(\ref{E:FD1})--(\ref{E:FD3}).

Let $\overline{x}$ be the class of $x\in \Gamma_{n+1} \Z(2)$ in $(G\Gamma
\Z(2))_n$. For $S^2$ we consider the simplicial model consisting in one generator
$y$ in degree two and only degeneracies above. The twisting morphism
$\eta:S^2\to G\Gamma \Z(2)$ is defined by the  relations~:
\begin{eqnarray*}
  \eta_{0}(*) &=& e,\qquad  \eta_{1}(*) = e\\
  \eta_{2}(y) &=& \bar{1} \\
  \eta_{3}(s_{1}y) &=& \overline{\sigma_{1}1}\\
  \eta_{3}(s_{2}y) &=& \overline{\sigma_{2}1}\\
  \eta_{3}(s_{0}y) &=& e \\
  \eta_{4}(s_{2}s_{1}y) &=& \overline{\sigma_{2}\sigma_{1}1}\\
  \eta_{4}(s_{3}s_{2}y) &=& \overline{\sigma_{3}\sigma_{2}1}\\
  \eta_{4}(s_{3}s_{1}y) &=& \overline{\sigma_{3}\sigma_{1}1}\\
  \eta_{4}(s_{i}s_{0}y) &=& e \qquad\text{for } 0\leq i<4 \\
   &\vdots&
\end{eqnarray*}
where $\bar{1}$ is a generator of 
$\mathcal{F}\left\{ \Z\backslash\{0\}\right\}$. In general, for
$n\geq 2$
$$
\eta_{n} (s_{j_{n-2}}\ldots s_{j_{1}}y)=\begin{cases}
e & \text{if }s_{j_{1}}=s_{0} \\
\overline{\sigma_{j_{n-2}}\ldots\sigma_{j_{1}}1} & \text{otherwise}
\end{cases}
$$
where $s_{j_{n-2}}\ldots s_{j_{1}}y$ is written in the canonical form.

The map $\eta$ is then determined by its value on the generator
$y$ of the model of
$S^2$, as is clear from the formula (\ref{E:FD3}).

\section{The simplicial model for the Hopf map}\label{S:main}
We now have all the tools to build our simplicial model for $S^3$. Denote by $\Z$
the set of integers and by
$\Z(2)$ the chain complex of abelian groups consisting in one copy of $\Z$ in
degree two and $0$ elsewhere. We apply the functor $\Gamma$ to get a
simplicial abelian group $\Gamma \Z(2)$. Therefore, by remark \ref{R:PTCP},
$$
p:G\Gamma \Z(2) \times_\eta S^2\to S^2
$$
is a principal twisted cartesian product whose fiber is $G\Gamma \Z(2)$ acting on
itself by left multiplication.   The map
$\eta: S^2\to G\Gamma \Z(2)$ is explained in the previous section.

\begin{theorem}
Let $S^2$ be endowed with the above simplicial model. A simplicial model for the Hopf
map
$S^3\to S^2$ is then given by the principal twisted cartesian product
$$
p:G\Gamma \Z(2) \times_\eta S^2\to S^2.
$$
\end{theorem}

\begin{proof}
The fibration $p:G\Gamma \Z(2) \times_\eta S^2\to S^2$ is a model for an element
of the set of $S^1$-bundles of base $S^2$, which contains the Hopf map. Now,
$S^1$-bundles of base $S^2$ are classified by $\Z$, and the Hopf map corresponds
to the class $1\in\Z$. All we have to show  is that our model $G\Gamma \Z(2)
\times_\eta S^2\to S^2$ corresponds indeed to the class $1\in \Z$. Consider the
diagramm
$$
\xymatrix{
		G\Gamma\Z(2) \ar[d]  & G\Gamma\Z(2) \ar[d]  &  \\
G\Gamma\Z(2)\times_\eta S^2 \ar[d]\ar[r]  & G\Gamma\Z(2)\times_\zeta \Gamma\Z(2)
\ar[d]  & \\
S^2 \ar[r]^\alpha & \Gamma\Z(2) \ar[r]^\beta & G\Gamma\Z(2) \\
		}
$$
where the two columns are fibrations and the composition $\beta\alpha$ is
the twisting map $\eta$. Recall from last section that the bottom composition
$\eta$ sends the generator $y$ of $S^2$ to the class of the generator $1\in\Z$.
Note that $ G\Gamma\Z(2)\times_\zeta \Gamma\Z(2)$ is acyclic and that the first
vertical fibration is classified by the map $\beta\alpha=\eta$. By choosing
$\eta$ to send a generator of $S^2$ to the generator $1\in\Z$ we guaranty that
our fibration $G\Gamma \Z(2) \times_\eta S^2\to S^2$ lies in the same class as the
Hopf map does and hence is a model of the later.
\end{proof}

\begin{remark}
In the previous proof, if we choose to send $y\in S^2$ to $m 1\in\Z$ via the map
$\alpha$, our fibration can model any $S^1$-bundle of base $S^2$ by letting $m$
vary over
$\Z$.

\end{remark}

\end{document}